\documentstyle[times,mssymb,12pt]{article}
\newcommand{\Vlad}{{Vl\u{a}du\c{t}}}
\newcommand{\DV}{{Drinfeld-$\!$\Vlad}}
\def\be{\begin{equation}}
\def\ee{\end{equation}}
\def\bea{\begin{eqnarray}}
\def\eea{\end{eqnarray}}
\def\Q{{\bf Q}}
\def\Z{{\bf Z}}
\def\fZ{{\footnotesize\bf Z}}

\def\C{{\bf C}}
\def\fC{{\footnotesize\bf C}}
\def\F{{\bf F}}
\def\HH{{\cal H}}
\def\PP{{\bf P}}
\def\X{{\rm X}}
\def\XX{{\cal X}}
\def\Xo{\X_0}
\def\XXo{\XX_0}
\def\0{^{\phantom0}}
\def\ra{\rightarrow}
\def\lra{\leftrightarrow}
\def\llra{\longleftrightarrow}
\def\geqs{\geqslant}

\begin{document}

\begin{center}

{\Large Explicit modular towers}

Noam D. Elkies\\
Department of Mathematics, Harvard University\\
Cambridge, MA 02138 (e-mail: {\tt elkies@math.harvard.edu})

\end{center}

\begin{quote}
{\small
{\bf Abstract.}
We give a general recipe for explicitly constructing asymptotically
optimal towers of modular curves such as $\{\Xo(l^n)\}_{n>1}$.
We illustrate the method by giving equations for eight towers
with various geometric features.  We conclude by observing that
such towers are all of a specific recursive form and speculate
that perhaps every tower of this form which attains the \DV\ bound
is modular.
}
\end{quote}

\vspace*{3ex}

{\bf Introduction.}  Explicit equations for modular curves
have attracted interest at least since the classical work of Fricke
and Klein.  Recent renewed interest in such equations has been
stimulated on the one hand by the availability of software for
symbolic computation and on the other hand by specific applications.
In~[E] we considered the use of modular curves to count rational
points on elliptic curves over large finite fields, and illustrated
some other applications of equations for the curves $\Xo(N)$ with
small~$N$\/ (say $N<10^3$).  Another kind of application is to
coding theory: good Goppa codes~[G] require curves of large genus~$g$
over a fixed finite field $k=\F_q$ whose number of rational points
grows as a positive multiple of~$g$.  Drinfeld and \Vlad\ showed that
as $g\ra\infty$ no multiple greater than $(q^{1/2}-1)g$ is possible.
Ihara~[I] and, independently, Tsfasman, \Vlad, and Zink~[TVZ] showed
that this upper bound is attained by the supersingular points on
appropriate modular curves when $q$ is a square.
For this application modular curves --- elliptic, Shimura, or Drinfeld
--- are needed whose level is too
high to apply the methods of~[E] directly, and in general one does
not expect to have any pleasant model for a curve of high genus.
However, if the curve is of smooth level then it tops a tower of
$O(\log q)$ covers of low degree, and one may hope to obtain
equations for the curve by writing those covers explicitly.

In this paper we show how to do this recursively for towers such as
$\{\Xo(l^n)\}_{n>1}$.  It turns out that only information about
the first few levels of the tower is needed, and that this information
can be obtained for modular elliptic curves using the methods of~[E],
and for some Shimura curves using only the ramification structure.
We then illustrate the method by giving
explicit formulas for eight asymptotically optimal towers:
six of elliptic modular curves, namely $\Xo(l^n)$ for $l=2,3,4,5,6$,
and $\Xo(3\cdot 2^n)$; and two of Shimura modular curves.
Over any finite field whose characteristic does not divide
the level of these modular curves, the towers are tamely ramified,
making it easy to calculate the genus of every curve in the tower.
[This contrasts with the wildly ramified tower of [GS1], whose
genus computation required some ingenuity; we show elsewhere
that that tower too is modular, of Drinfeld type.]
For each finite field~$k$\/ over which one of our towers is
asymptotically optimal, the optimality can then be shown by
elementary means, independent of the tower's modular provenance,
by exhibiting the coordinates of the rational (supersingular)
points.  These formulas may also have other uses, e.g.\ in finding
explicit modular parametrizations of elliptic curves with smooth
conductor, or in connection with generalizations of the
arithmetic-geometric mean (which corresponds to the $\Xo(2^n)$ tower)
as in~[S1,S2]; we hope to pursue these connections in future papers.
We conclude this paper with a speculation
concerning the modularity of ``any'' asymptotically optimal tower.

{\bf The curves $\Xo(l^n)$.}  Fix a prime $l>1$.  For positive $n$,
the elliptic modular curve $\Xo(l^n)$ over any field~$k$\/ in which
$l\neq0$ parametrizes elliptic curves
with a cyclic $l^n$-isogeny, or equivalently sequences of
$l$-isogenies
\be
E_0 \ra E_1 \ra E_2 \ra \cdots \ra E_n
\label{lisogs}
\ee
such that the composite isogeny $E_{j-1} \ra E_{j+1}$ of degree $l^2$
is cyclic for each $j$ with $0<j<n$.  Thus for each $m=0,1,\ldots,n$
there are $n+1-m$ maps $\pi_j: \Xo(l^n)\ra\Xo(l^m)$ obtained by
extracting for some $j=0,1,\ldots,n-m$ the cyclic $l^m$-isogeny
$E_j \ra E_{j+m}$ from (\ref{lisogs}).  Each of these maps has degree
$l^{n-m}$, unless $m=0$ when the degree is $(l+1) l^{n-1}$.
In particular we have a tower of maps
\def\pora{\stackrel{\pi_0}{\ra}}
\be
\Xo(l^n) \pora \Xo(l^{n-1}) \pora \Xo(l^{n-2}) \pora \cdots
\pora \Xo(l^2) \pora \Xo(l),
\label{Xotower}
\ee
each map being of degree~$l$.  Each $\Xo(l^n)$ also has an
Atkin-Lehner involution $w_l=w_l^{(n)}$, taking a cyclic
$l^n$-isogeny to its dual isogeny, and the sequence (\ref{lisogs})
to the sequence
\be
E_n \ra \cdots \ra E_2 \ra E_1 \ra E_0
\label{sgosil}
\ee
of dual isogenies.  We thus have
\be
w_l^{(m)} \circ \pi_j = \pi_{n-m-j} \circ w_l^{(n)},
\label{pi,w}
\ee
where $\pi_j, \pi_{n-m-j}$ are our $j$th and $(n-m-j)$th
maps from $\Xo(l^n)$ to $\Xo(l^m)$.

When $k=\C$, we may regard $\Xo(N)$ as the quotient of the extended
upper half-plane $\HH^*= \HH \cup \PP^1(\Q)$ by $\Gamma_0(N)$.
Then the $\Gamma_0(l^n)$ orbit of a point $\tau\in\HH$ parametrizes
the isogeny between the elliptic curves with period lattices
$\Z + \tau \Z$ and $l^{-n} \Z + \tau \Z$,
the map $\pi_j$ takes the $\Gamma_0(l^n)$ orbit
of~$\tau$ to the $\Gamma_0(l^m)$ orbit of $l^j \tau$, and
the involution $w_l^{(n)}$ is represented by $\tau \llra -1/l^n \tau$.

Now the key observation is that explicit formulas for
$\Xo(l),\Xo(l^2)$, together with the involutions
$w_l^{(1)},w_l^{(2)}$ of these curves and the map
$\pi_0:\Xo(l^2)\ra\Xo(l)$ between them, suffice to exhibit
the entire tower (\ref{Xotower}) explicitly:

{\bf Proposition.}  {\sl For $n\geqs2$ the product map
\be
\pi = \pi_0 \times \pi_1 \times \pi_2 \times \cdots \times \pi_{n-2}:
\Xo(l^n) \ra (\Xo(l^2))^{n-1}
\label{prodmap}
\ee
is a 1:1 map from $\Xo(l^n)$ to the set of
$(P_1,P_2,\ldots,P_{n-1}) \in \bigl(\Xo(l^2)\bigr)^{n-1}$ such that
\be
\pi_0 \bigl( w_l^{(2)} (P_j) \bigr) =
w_l^{(1)} \bigl( \pi_0(P_{j+1}) \bigr)
\label{conds}
\ee
for each $j=1,2,\ldots,n-2$.}

Informally speaking, we get from $\Xo(l^2)$ up to $\Xo(l^n)$ by
iterating $n-2$ times the involution $w_l^{(2)}$ composed with the
``$l$-valued involution'' $\pi_0^{-1} w_l^{(1)} \pi_0\0$.
Of course the maps $\pi_j: \Xo(l^n) \ra \Xo(l^m)$ (for $m\geqs 2$)
are then simply
\be
(P_1,\ldots,P_{n-1}) \mapsto (P_{j+1},\ldots,P_{j+m-1}),
\label{pij}
\ee
and the involution $w_l^{(n)}$ is
\be
(P_1,P_2,\ldots,P_{n-2},P_{n-1}) \llra
(w_l^{(2)}P_{n-1},w_l^{(2)}P_{n-2},\ldots,w_l^{(2)}P_2,w_l^{(2)}P_1),
\label{w}
\ee
i.e.\ reversing the order of $P_1,\ldots,P_{n-1}$ and applying
$w_l^{(2)}$ to each coordinate.

{\sl Proof}\/\footnote{
  More properly, a proof sketch, since we suppress some details,
  such as what happens at the cusps.  To show that our formulas
  extend to the cusps one may either quote general facts about
  maps between affine and projective algebraic curves, or regard
  the cusps as parametrizing isogenies between Tate curves.
  Also, two cyclic $l$-isogenies may determine the same point
  on $\Xo(l)$ without being isomorphic; the necessary and sufficient
  condition is that they become isomorphic over the algebraic closure.
  }:
That the map is 1:1 to its image is clear, because
a sequence (\ref{lisogs}) of $l$-isogenies is determined by the
$l^2$-isogenies $E_{j-1} \ra E_{j+1}$ parametrized by the
$j$th coordinate of~$\pi$ ($0<j<n$).  Now $(P_1,\ldots,P_{n-1})$
is in the image of~$\pi$ if and only if the $l^2$-isogenies
parametrized by $P_1,\ldots,P_{n-1}$, regarded as sequences
$E_0^j \ra E_1^j \ra E_2^j$ of $l$-isogenies, fit together
to form a sequence (\ref{lisogs}) with $E_i^j=E_{i+j}$,
i.e.\ if and only if the isogenies $E_1^j \ra E_2^j$ and
$E_0^{j+1} \ra E_1^{j+1}$ coincide for each $j=1,2,\ldots,n-2$.
But these isogenies are represented by the points $\pi_1(P_j)$
and $\pi_0(P_{j+1})$ on $\Xo(l)$.
Thus the necessary and sufficient condition is that
\be
\pi_1(P_j) = \pi_0(P_{j+1})
\label{cond1}
\ee
for each $j=1,2,\ldots,n-2$; applying $w_l^{(1)}$ to both sides,
and then (\ref{pi,w}) to $w_l^{(1)} \bigl( \pi_1(P_j) \bigr)$,
then yields the equivalent form (\ref{conds}).~~$\Box$

{\bf Examples: The cases $l=2,3,5$.}  Our formulas are
particularly simple when $\Xo(l^2)$ (and thus also $\Xo(l)$) has
genus~0, for then we may use a Hauptmodul (or for that matter any
rational parameter\footnote{A ``Hauptmodul'' is a rational parameter
with a pole of leading coefficient~1 at the infinite cusp, i.e.\
a degree-1 rational function of the form $q^{-1} + O(1)$.})
of $\Xo(l^2)$ to regard $P_1,\ldots,P_{n-1}$ as $n-1$
rational coordinates on~$\Xo(l^n)$, and
(\ref{conds}) as the $n-2$ algebraic relations on those coordinates
that determine the curve $\Xo(l^n)$.  This happens for $l=2,3,5$;
we exhibit formulas for each of these cases.

In the first two cases the cover $\pi_0: \Xo(l^2) \ra \Xo(l)$
is cyclic.\footnote{
  For any $N$, the cover $\pi_0: \Xo(N^2) \ra \Xo(N)$ is cyclic
  if and only if the unit group of \fZ$/N$\fZ\ has exponent~$2$,
  which happens when $N|24$.  When $N=1,2,3,4,6$, it is furthermore
  true that $\pm1$ are the only units of \fZ$/N$\fZ, and then
  the cyclic $N$-isogenies $E_1\ra E_0,E_2$ together with the
  Weil pairing on $E_1[N]$ determine a complete level-$N$\/
  structure on~$E_1$ mod~$\pm1$, i.e.\ $\Xo(l^2)\cong\X(l)$.
  }
For $l=2$ we parametrize $\Xo(l^2)=\Xo(4)$ by
\be
\xi(\tau) := 1 + \frac18
\left(\frac{\eta(\tau)}{\eta(4\tau)}\right)^{\!8}
= \frac18(q^{-1} + 20q - 62q^3 + 216q^5 - 641q^7 + - \cdots),
\label{xi4}
\ee
where as usual $q=e^{2\pi i \tau}$ and $\eta$ is the weight-$\frac12$
modular form $\prod_{r=1}^\infty (1-q^r)$.  Using the functional
equation
\be
\eta(-1/\tau) = (\tau/i)^{1/2} \eta(\tau),
\label{etafunceq}
\ee
we find that the involution $w_2^{(2)}$ takes $\xi(\tau)$ to
\be
\xi(-1/4\tau)
= 1 + 32 \left(\frac{\eta(4\tau)}{\eta(\tau)}\right)^{\!8}
= 1 + \frac 4 {\xi(\tau) - 1}
= \frac {\xi(\tau) + 3} {\xi(\tau) - 1} \, .
\label{w4}
\ee
Let $h_2$ be the $\Xo(2)$ Hauptmodul
\be
h_2(\tau) = \left(\frac{\eta(\tau)}{\eta(2\tau)}\right)^{\!24}
= q^{-1} - 24 + 276 q - 2048 q^2 + 11202 q^3 - + \cdots.
\label{h2}
\ee
Computing the map $\pi_0: \Xo(4) \ra \Xo(2)$ then amounts to
writing $h_2$ as a rational function in~$\xi$.  We do this
by in effect expanding this function as a continued fraction.
Necessarily that function has degree~2 with a simple pole at 
the cusp $\xi=\infty$.  But then $h_2 - 8\xi + 24$ is a
rational function of degree~1 in~$\xi$ with a simple zero
at~$\infty$, i.e.\ the inverse of a polynomial of degree~1.
Comparing the $q$-expansions of $1/(h_2 - 8\xi + 24)$ and~$\xi$,
we find that this polynomial is $(\xi+1)/32$ and recover the
formula
\be
h_2(\tau) = 8 \frac{ (\xi(\tau)+1)^2 } { (\xi(\tau)-1) }.
\label{pi0:42}
\ee
Using our formula (\ref{w4}) for the involution $w_2^{(2)}$
we then obtain also
\be
h_2(\tau) = \frac{64}{ \xi(-1/4\tau)^2 - 1}.
\label{pi1:42}
\ee
But $w_2^{(1)}$ acts on $\Xo(2)$ by $h_2 \lra 2^{12} / h_2$
(again by (\ref{etafunceq})).  Thus $h_2(2\tau)$ is both
$64(\xi(\tau)^2 - 1)$ and $64/(\xi(-1/8\tau)^2 - 1)$.
Equating these two expressions yields an equation for
the modular curve $\Xo(8)$; more generally we now deduce
from our Proposition the following explicit equations for the
modular curve $\Xo(2^n)$ for each $n>1$:

Let $x_j$ ($0<j<n$) be the rational function $\xi(2^{j-1}\tau)$
on that curve (this is the coordinate $P_j$ of the Proposition);
then $(x_1,\ldots,x_{n-1})$ identifies $\Xo(2^n)$ with the
curve in $(\PP^1)^{n-1}$ specified by the $n-2$ equations
\be
(x_j^2-1) (z_{j+1}^2-1) = 1
\qquad (j=1,\ldots,n-2),
\label{xz2}
\ee
where
\be
z_j := (x_j+3)/(x_j-1)
\label{z2}
\ee
is obtained from $x_j$ by the involution $w_2^{(2)}$.

Curiously we obtain analogous equations for $\Xo(3^n)$ by
replacing the exponent~2 by~3 in~(\ref{xz2}) and, as if
to compensate, changing the constant term~3 to~2 in~(\ref{z2}):
the curve $\Xo(3^n)$ is isomorphic with the locus of
$(x_1,\ldots,x_{n-1})$ in $(\PP^1)^{n-1}$ satisfying
\be
(x_j^3-1) (z_{j+1}^3-1) = 1
\qquad (j=1,\ldots,n-2),
\label{xz3}
\ee
where
\be
z_j := (x_j+2)/(x_j-1).
\label{z3}
\ee
Here the coordinate functions $x_j$ on $\Xo(3^n)$ are
$\xi(3^{j-1}\tau)$, where
\be
\xi(\tau) =
1 + \frac13 \left(\frac{\eta(\tau)}{\eta(9\tau)}\right)^{\!3}
= \frac13(q^{-1} + 5q - 7q^5 + 3q^8 + 15q^{11} - 32q^{14}
\cdots),
\label{xi9}
\ee
so $\xi$ generates the field of rational functions on $\Xo(9)$.
The involution $w_3^{(2)}$ takes this $\xi$ to
$1 + 3/(\xi-1) = (\xi+2)/(\xi-1)$, whence (\ref{z3});
the Hauptmodul
\be
h_3(\tau) = \left(\frac{\eta(\tau)}{\eta(3\tau)}\right)^{\!12}
= q^{-1} - 12 + 54 q - 76 q^2 - 243 q^3 + 1188 q^4 \cdots
\label{h3}
\ee
goes to $3^6/h_3$ under $w_3^{(1)}$, and $h_3(3\tau)$ is both
$27/(\xi(-1/27\tau)^3 - 1)$ and $27(\xi(\tau)^3 - 1)$, from
which the equations (\ref{xz3},\ref{z3}) for $\Xo(3^n)$ follow
thanks to our Proposition.

Note that the already simple equations for $\Xo(2^k)$, $\Xo(3^k)$
simplify even further in characteristic 3,~2 respectively:
taking $y_j = 1 - x_j^{-1}$ in (\ref{xz2},\ref{z2}) and setting
$3=0$ yields
\be
y_{j+1}^2 = y_j\0-y_j^2 ,
\label{2mod3}
\ee
and the same substitution in (\ref{xz3},\ref{z3}) with $2=0$ produces
\be
y_{j+1}^3 = y_j^3 + y_j^2 + y_j\0.
\label{3mod2}
\ee
In this guise these asymptotically optimal towers were obtained by
Garcia and Stichtenoth~[GS2, Examples C,D], independent (as in [GS1])
of their modular interpretation.  In both cases the supersingular
points are the poles of $y_1$ and thus of all the $y_j$.

Finally for $l=5$ we obtain
\be
P(x_j) P(z_{j+1}) = 125
\qquad (j=1,\ldots,n-2),
\label{xz5}
\ee
where
\be
P(X):= X^5 + 5 X^3 + 5 X - 11,
\qquad
z_j := (x_j+4)/(x_j-1).
\label{z5}
\ee
Here the coordinate $x_j$ is $\xi(5^{j-1}\tau)$ where
\be
\xi(\tau) =
1 + \frac{\eta(\tau)}{\eta(25\tau)} =
q^{-1} - q + q^4 + q^6 - q^{11} - q^{14} + q^{21} + q^{24} - q^{26}
\cdots.
\label{xi25}
\ee
As usual, $z_j$ is the image of $x_j$ under $w_l^{(2)}$,
and $P(x_j)=h_5(5\tau)$ where $h_5$ is the $\Xo(5)$ Hauptmodul
\be
h_5(\tau) = \left(\frac{\eta(\tau)}{\eta(5\tau)}\right)^{\!6}
= q^{-1} - 6 + 9 q + 10 q^2 - 30 q^3 + 6 q^4 - 25 q^5 \cdots
\label{h5}
\ee
with $h_5(\tau) h_5(-1/5\tau) = 125$.  The polynomial $P(X)$ is
necessarily not as simple as the polynomials $X^2-1$, $X^3-1$
occurring in (\ref{xz2},\ref{xz3}), because the cover
$\Xo(25)\ra\Xo(5)$ is not cyclic.  It is, however, dihedral,
as may be seen from the fact that $P(W-1/W)=W^5-11-W^{-5}$.

{\bf First variation: composite $l$.}  The assumption that
$l$\/ be prime was not necessary; the entire description
carries over to the composite case, except for the incidental
point that the degree of the maps $\pi_j: \Xo(l^n) \ra \X(1)$
is given by a formula more complicated than $(l+1)l^{n-1}$
[namely $l^n \prod_{p|l} (1+\frac1p)$].  For instance we
exhibit formulas for the cases $l=4,6$, where the cover
$\pi_0: \Xo(l^2) \ra \Xo(l)$ is still cyclic.

In the first case $l=4$ the curve $\Xo(l^2)=\Xo(16)$ is still
rational, and we obtain formulas remarkably similar to those
for $l=2,3$ by choosing
\be
\xi(\tau) = 1 + \frac12
\frac{\eta^2(\tau)\eta(8\tau)}{\eta(2\tau)\eta^2(16\tau)}
= \frac12(
q^{-1} + 2 q^3 - q^7 - 2 q^{11} + 3 q^{15} + 2 q^{19}
\cdots )
\label{xi16}
\ee
as a rational coordinate on $\Xo(16)$.  Then $w_4^{(2)}$ takes $\xi$ to
$(\xi+1)/(\xi-1)$.  The $\Xo(4)$ Hauptmodul
\be
h_4(\tau) = \left(\frac{\eta(\tau)}{\eta(4\tau)}\right)^{\!8}
= q^{-1} - 8 + 20 q - 62 q^3 + 216 q^5 - 641 q^7 + - \cdots
\label{h4}
\ee
(cf.~(\ref{xi4})) is mapped by $w_4^{(1)}$ to $4^4/h_4$, and we compute
\be
h_4(4\tau) = \frac{16} {\xi(-1/64\tau)^4 - 1} 
= 16 \bigl(\xi(\tau)^4 - 1\bigr).
\label{16to4}
\ee
Therefore $\Xo(4^n)$ is isomorphic with the locus of
$(x_1,\ldots,x_{n-1})$ in $(\PP^1)^{n-1}$ satisfying
\be
(x_j^4-1) (z_{j+1}^4-1) = 1
\qquad (j=1,\ldots,n-2),
\label{xz4}
\ee
where
\be
z_j := (x_j+1)/(x_j-1),
\label{z4}
\ee
the coordinate functions $x_j$ on $\Xo(4^n)$ being $\xi(4^{j-1}\tau)$.
Of course the resulting curves also occur in the $\Xo(2^n)$ tower,
but this fact is far from obvious from comparison of the formulas
(\ref{xz2},\ref{z2}) and (\ref{xz4},\ref{z4}).

The case of $l=6$ is slightly more complicated because
the curve $\Xo(l^2)$ is no longer rational.  It is, however,
an elliptic curve with a simple Weierstrass equation: the
ring of rational functions on $\Xo(36)$ regular except
possibly at the cusp $\tau=i\infty$ is generated by
\bea
\xi(\tau) &=&
\frac{\eta(12\tau)\eta^3(18\tau)}{\eta(6\tau)\eta^3(36\tau)}
\; = \; q^{-2} + q^2 + q^8 - q^{14} - q^{20} + q^{26} + 2 q^{32} \cdots,
\label{xi36} \\
\gamma(\tau) &=&
\frac{\eta^4(12\tau)\eta^2(18\tau)}{\eta^2(6\tau)\eta^4(36\tau)}
\; = \; q^{-3} + 2 q^3 + q^9 - 2 q^{15} - 2 q^{21} + 2 q^{27} + 4 q^{33}
\cdots,
\label{gamma36}
\eea
related by the Weierstrass equation
\be
\gamma^2 = \xi^3 + 1.
\label{Xo36}
\ee
The involution $w_6^{(2)}$ has a fixed point at $i/6$.
An involution of an elliptic curve which has a fixed point
must be multiplication by~$-1$ composed with a translation.
Thus to determine $w_6^{(2)}$ we need only find the image of
one point.  It is easiest to do this with the cusp $\tau=i\infty$:
its image is the cusp $\tau=0$, at which $(\xi,\gamma)=(2,3)$
(a 6-torsion point on the curve (\ref{Xo36})).  It remains only
to find the map from $\Xo(36)$ to $\Xo(6)$ and the involution
$w_6^{(1)}$.  We use the Hauptmodul
\be
h_6(\tau) = \frac{\eta^5(\tau)\eta(3\tau)}{\eta(2\tau)\eta^5(6\tau)}
= q^{-1} - 5 + 6 q + 4 q^2 - 3 q^3 - 12 q^4 - 8 q^5 + 12 q^6 \cdots.
\label{h6}
\ee
Then $w_6^{(1)}$ takes $h_6$ to $72/h_6$, and by comparing
$q$-expansions we find
\be
h_6(6\tau) = \xi^3(\tau) - 8.
\label{36to6}
\ee
We thus identify $\Xo(6^n)$ with the curve of $(n-1)$-tuples
$\bigl( (x_1,y_1), \ldots, (x_{n-1},y_{n-1}) \bigr)$
of points on the elliptic curve $y^2=x^3+1$ satisfying
the $n-2$ conditions
\be
(x_j^3-8) (z_{j+1}^3-8) = 72
\qquad (j=1,\ldots,n-2),
\label{xz6}
\ee
where
\be
z_j := \left( \frac{y_j+3}{x_j-2} \right)^{\!2} - x_j - 2
\label{z6}
\ee
is the $x$-coordinate of the point $(2,3) - (x,y)$ on $y^2=x^3+1$.
Unlike the curves in the $\Xo(4^n)$ tower, these curves $\Xo(6^n)$ 
are new; to be sure they could also be exhibited as composita of
the already known covers $\Xo(2^n)/\X(1)$ and $\Xo(3^n)/\X(1)$,
but those models are much harder to work with because of the
complicated singularities above the branch points $j=0,12^3,\infty$.

{\bf Second variation: changing the base of the tower.}
Instead of the tower of modular curves $\Xo(l^n) = \HH^*/\Gamma_0(l^n)$
we could use $\HH^*/(\Delta \cap \Gamma_0(l^n))$ where $\Delta$ is
some other congruence subgroup of PGL$_2(\Q)$, as long as the
modulus of the congruence is prime to~$l$.  For instance, given $N>1$
with $(l,N)=1$ we could use the tower $\Xo(N l^n)$ of curves
parametrizing sequences of \hbox{$l$-isogenies} between pairs of elliptic
curves related by a cyclic $N$-isogeny.  Again these curves with $n>1$
form a tower related by maps $\pi_j$ of $l$-power degree and admitting
involutions $w_l^{(n)}$, and knowing these maps and involutions for
$n=1,2$ yields explicit formulas for $\Xo(N l^n)$ for all $n>1$ as
in our Proposition.\footnote{
  As with $\Xo(6^n)$ we could also obtain $\Xo(Nl^n)$ as a compositum
  of $\Xo(N)$ and $\Xo(l^n)$, but the resulting model is highly
  singular.  Warning: on $\Xo(Nl^n)/$\fC\ the involution
  $w_l^{(n)}$ is given not by $\tau\lra-1/l^n\tau$ but by a fractional
  linear transformation of the same determinant that reduces mod~$N$\/
  to an element of $\Gamma_0(N)$.  We do still have a simple formula
  $\tau\lra-1/N l^n\tau$ for the product of $w_l^{(n)}$ with the
  Atkin-Lehner involution $w_N$.
  }

We illustrate with the case $l=2$, $N=3$.  In this case the first
two curves $\Xo(6),\Xo(12)$ are rational and we can mimic our
procedure for the towers $\Xo(l^n)$ with $l=2,3,4$.  Our $(n-1)$
coordinates on $\Xo(3\cdot 2^n)$ will be $x_j = \xi(2^{j-1}\tau)$
($0<j<n$) where
\be
\xi(\tau) =
\frac{\eta^4(4\tau)\eta^2(6\tau)}{\eta^2(2\tau)\eta^4(12\tau)}
= q^{-1} + 2 q + q^3 - 2 q^5 - 2 q^7 + 2 q^9 + 4 q^{11} \cdots
\label{xi12}
\ee
(cf.~(\ref{gamma36})) is a Hauptmodul for $\Xo(12)$.  It is this
time more convenient to let $h_6$ be the $\Xo(6)$ Hauptmodul
\be
h_6(\tau) = \left(
\frac{\eta(2\tau)\eta^3(3\tau)}{\eta(\tau)\eta^3(6\tau)}
\right)^{\!3}
= q^{-1} + 3 + 6 q + 4 q^2 - 3 q^3 - 12 q^4 - 8 q^5
\cdots,
\label{h6'}
\ee
which differs by~$8$ from our choice in~(\ref{h6}).
We may then represent $w_2^{(1)}$ and $w_2^{(2)}$ by
$\tau\llra(2\tau-1)/(6\tau-2)$ and $\tau\llra(4\tau+3)/(4\tau+4)$;
these involutions take $h_6$ to $-8/h_6$ and $\xi$ to $(3-\xi)/(1+\xi)$.
By computing the quadratic map $\Xo(12)\ra\Xo(6)$ we find that this time
\be
h_6(2\tau) = \frac{-8} {\xi(w_2^{(2)}\tau)^2 - 1} 
= \xi(\tau)^2 - 1.
\label{12to6}
\ee
Thus the equations on $x_1,\ldots,x_{n-1}$ defining $\Xo(3\cdot 2^n)$
are
\be
(x_j^2-1) (z_{j+1}^2-1) = -8
\qquad (j=1,\ldots,n-2),
\label{xz2.3}
\ee
where
\be
z_j := (3-x_j)/(1+x_j).
\label{z2.3}
\ee
Note that in this case the curves in our tower also have
an involution $w_3$ commuting with all the $w_2^{(n)}$;
we find that this involution is $x_j\lra -3/x_j$.  That
this in fact acts on our model of $\Xo(3\cdot 2^n)$ is
easy to check after writing (\ref{xz2.3},\ref{z2.3}) in
the equivalent form
\be
(x_{j+1}-1) x_j^2 = x_{j+1}^2 + 3 x_{j+1}.
\label{xx2.3}
\ee

{\bf Third variation: Shimura modular curves.}  Shimura curves
generalize the classical elliptic modular curves: instead of
$\HH^*/\Gamma$ for an arithmetic subgroup~$\Gamma$ of PGL$_2(\Q)$,
they are the quotients $\HH/\Gamma$ by an arithmetic subgroup of
a quaternion algebra~$A$\/ over some totally real number field~$K$,
with $A$ ramified at all but one of the infinite places of~$K$.
Instead of elliptic curves, these Shimura curves parametrize
principally polarized abelian varieties with endomorphisms
by~$A$ and some extra structure determined by the choice of~$\Gamma$.
There are Shimura curves $\XXo(I)$ ($I$\/ an ideal of~$K$\/ coprime
with the discriminant of~$A$) analogous to $\Xo(N)$, which have
Atkin-Lehner involutions and form towers, and whose reductions
at a prime of~$K$\/ are asymptotically optimal over the quadratic
extensions of its residue field.  These towers may be obtained
from their first two levels by the recipe of our Propositions.

Unlike the classical $\Xo(N)$, the analogous Shimura curves $\XXo(I)$
have no cusps.  Thus the curves and maps between them cannot be
computed using $q$-expansions.  Even worse, in general we do not even
have explicit equations for the abelian varieties parametrized by these
curves.  Nevertheless we can in many cases use the ramification
behavior of the covers to determine the necessary maps completely.
We illustrate this with two examples which have the additional
feature of involving only cyclic covers which become unramified after
finitely many steps and thus also occur in class-field towers.\footnote{
  Note that this is not possible with elliptic modular curves,
  or for that matter with Drinfeld modular curves, precisely because
  of their cusps.  However, the ramification in towers of elliptic
  or Drinfeld modular curves is small enough to be captured by
  a tower of ray class fields, suggesting that ray class-field
  towers might be a fruitful source of curves with many points
  even over finite fields of non-square order.}

We start with $K,A$ such that $A^*$ contains an arithmetic subgroup
$\Delta$ which is also a triangle group.  Such $\Delta$ have been
classified completely~[T]: there are 76, in 18 quaternion algebras
(not including the nine triangle subgroups of PGL$_2(\Q)$ with
one or more cusps among the vertices).  For our first example we
take $K=\Q(\sqrt{3}\,)$ and $A/K=$ the quaternion algebra ramified at
$(\sqrt3\,)$ and at one infinite place, and choose for $\Delta$
the group called $\Gamma^{(+)}(A,O_1) = \Gamma^{(*)}(A,O_1)$ in~[T],
which is identified there with the $(2,4,12)$ triangle group.
We shall construct the tower $\{\XXo(\wp_2^n)\}_{n>1}$, where
$\wp_2$ is the prime of~$K$\/ of residue field~$\F_2$.

The curve $\XX(1)=\HH/\Delta$ is rational.  We choose a
coordinate~$J$\/ taking the values $1,0,\infty$ at the
elliptic points of order $2,4,12$.  The curve $\XXo(\wp_2)$
consists of ordered pairs of points of $\XX(1)$ related by a
``$\wp_2$-isogeny''; choosing one of these points yields the
degree-3 map $\pi_0: \XXo(\wp_2) \ra \XX(1)$.  We next
determine the ramification of this map.  In general,
the map $\XXo(I) \ra \XX(1)$ is branched only above
elliptic points of $\XX(1)$, if a point~$P$\/ of $\XXo(I)$
above an elliptic point of order~$e$ parametrizes an
isogeny to some other point of order $e'$ then the ramification
index at~$P$\/ is the denominator of the fraction $e'/e$.
[A non-elliptic point is taken to have order~1.]
We may regard $\XXo(\wp_2)$ as a symmetric $(3,3)$ correspondence
on $\XX(1)\times\XX(1)$.  We then see that the point $J=\infty$
of order~12 must correspond to the point $J=0$ of order~4
with multiplicity~3; the point $J=0$ corresponds to $J=1$
doubly and $J=\infty$ singly; and $J=1$ corresponds to
$J=0$ singly and some other point doubly.  No other points
of $\XX(1)$ are ramified in $\XXo(\wp_2)$.  Thus by the 
Riemann-Hurwitz formula $\XXo(\wp_2)$ is again a rational
curve, and $J$\/ is a function of degree~3 with a triple
pole such that $J$\/ and $J-1$ both have double zeros.
Up to Aut$(\PP^1)$ there is a unique such function;
we choose a rational coordinate~$t$\/ on~$\XXo(\wp_2)$
such that $J=t(4t-3)^2$ (so $J-1=(t-1)(4t-1)^2$).
Then the involution\footnote{
  We suppress the unwieldy subscript $\wp_2$.
  }
$w^{(1)}$ must interchange the points $t=0$, $t=\infty$
parametrizing isogenies between $J=0$ and $J=\infty$,
and the points $t=1$, $t=3/4$ parametrizing isogenies
between $J=0$ and $J=1$.  Therefore $w^{(1)}(t)=3/4t$.

Now the curve $\XXo(\wp_2^2)$ covers $\XXo(\wp_2)$ with
degree~2, and the only branch points are $t=\infty$ and
$t=3/4$.  Thus $\XXo(\wp_2^2)$ is again a rational
curve, and we may choose a rational coordinate~$\xi$
for it such that $t=(\xi^2+3)/4$.  Of the points $\xi=\pm1$
above $t=1$, one must parametrize a $\wp_2^2$-isogeny from $J=1$
to $J=\infty$, the other an isogeny from $J=1$ to itself; we
choose $\xi$ so the former point is $\xi=1$.  Then $w^{(2)}$
must switch that point with the point $\xi=\infty$, and
fix the other point $\xi=-1$; therefore this involution
is $\xi \llra (\xi+3)/(\xi-1)$.  As a further check on
the computation, note that this involution also switches
the two points $\xi=\pm\sqrt{-3}$ above $t=0$, parametrizing a
$\wp_2^2$-isogeny from $J=0$ to itself.

We now have all the information needed to determine the
Shimura modular curve $\XXo(\wp_2^n)$ for all $n>1$:
that curve has $n-1$ coordinates $x_1,\ldots,x_{n-1}$,
satisfying the $n-2$ relations
\be
\left( \frac{x_j^2+3}{4} \right)
\left( \frac{z_{j+1}^2+3}{4} \right)
= \frac34,
\label{s2tmp}
\ee
that is,
\be
(x_j^2+3) (z_{j+1}^2+3) = 12
\qquad (j=1,\ldots,n-2),
\label{shim:xz2}
\ee
where
\be
z_j := (x_j+3)/(x_j-1),
\label{shim:z2}
\ee
the same involution we used in (\ref{z2}) for the tower of
classical modular curves $\Xo(2^n)$.  Unlike these curves,
though, the Shimura tower $\XXo(\wp_2^n)$ turns out to
be unramified past $n=5$, as may be seen either directly
from the formulas (\ref{shim:xz2},\ref{shim:z2}) or from the
general description of ramification in the map
$\XXo(I) \ra \XX(1)$.  Since each step in the tower
is a cyclic extension, it follows that over any finite
field of odd characteristic the tower is dominated by
the 2-class-field tower of the curve $\XXo(\wp_2^5)$.

For our second example, we choose for $K$\/ the cubic field
$\Q(2\cos \pi/9)$ and for $A$ the quaternion algebra ramified
only at two of the three infinite places of~$K$.  Then we find
in~[T] that the group of units of norm~1 in~$A$ is the $(2,3,9)$
triangle group.  We exhibit the tower $\{\XXo(\wp_3^n)\}_{n>1}$,
where $\wp_3$ is the prime of~$K$\/ of residue field~$\F_3$.
The equations were obtained in the same way that we found
(\ref{shim:xz2},\ref{shim:z2}); we leave the intermediate
steps as an exercise.  Again we find formulas similar to
those we obtained earlier (\ref{xz3},\ref{z3}) for the
classical modular curves: there are $n-1$ coordinates
$x_1,\ldots,x_{n-1}$, related by $n-2$ equations
\be
x_j^3 + z_{j+1}^3 = 1
\qquad (j=1,\ldots,n-2)
\label{shim:xz3}
\ee
(this time even simpler than the equation (\ref{xz3}) for
the classical case), where again
\be
z_j := (x_j+2)/(x_j-1).
\label{shim:z3}
\ee
Again the tower has cyclic steps and is unramified after
finitely many steps; we find that it is dominated by
the 3-class-field tower of the curve $\XXo(\wp_3^4)$.

{\bf Fantasia: a speculation on modularity.}  All our towers are
of the following form: the bottom curve $C_1$ over some finite
field~$k$\/ is equipped with an irreducible correspondence
$\Phi\subset C_1 \times C_1$ of bidegree $(l,l)$ and a set
$S\subset C_1(k)$ of rational points each of which corresponds
under~$\Phi$ with $l$\/ {\sl distinct} points also in~$S$\/;
the $n$-th curve in the tower is then the curve $C_n$ of $n$-tuples
$(P_1,\ldots,P_n)\in C_1^n$ such that $(P_j,P_{j+1})\in\Phi$
for $j=1,2,\ldots,n-1$.  Then $C_n$ has at least $l^{n-1}|S|$
rational points, and at least when $\Phi$ is tamely ramified
we can find the genus of $C_n$ as a function of~$n$.  For instance,
if $C_1=\Xo(l^2)$, $\Phi$ is the image of $\Xo(l^3)$ under
$\pi_0 \times \pi_1$, and $S$\/ is the set of supersingular
points, then we recover the tower of curves $\Xo(l^{n+1})$
of our Proposition.

But the $(C_1,\Phi,S)$ description makes no assumption of modularity:
we can, as in~[GS2], try any $C_1$ and $\Phi$ and hope to find an $S$\/
that yields many points on $C_n$.  In fact, several such $(C_1,\Phi)$
were found to admit $S$\/ large enough to make the tower $\{C_n\}$
asymptotically optimal~[GS1,GS2].  However, in each such case $\{C_n\}$
was subsequently explained as a modular tower.

This leads us to speculate: perhaps every asymptotically optimal
tower of this recursive form must be modular?

{\bf Acknowledgements.} Thanks to H. Stichtenoth and P. Sol\'e
for e-mail correspondence on the explicit towers and AGM
connections eventually published in [GS2] and~[S1] respectively.
Thanks also to J. Cazaran for reading the manuscript and spotting
several local errors as the paper was going into print.
Symbolic computations were greatly facilitated by the computer
packages {\sc pari} and {\sc macsyma}.

This work was made possible in part by funding from the
National Science Foundation and the Packard Foundation.

{\large\bf References}

[E] Elkies, N.D.: Elliptic and modular curves over finite fields
and related computational issues.  To appear in the proceedings
of the conference {\em Computational Perspectives on Number Theory}
held September, 1995 in Chicago in honor of A.O.L. Atkin,
J. Teitelbaum, ed.

[G] Goppa, V.D.: Codes on algebraic curves, {\em Soviet Math.\ Dokl.}\
{\bf 24} (1981) \#1, 170--172.

[GS1]  Garcia, A., Stichtenoth, H.: A tower of Artin-Schreier
extensions of function fields attaining the Drinfeld-$\!$Vladut bound.
{\em Invent.\ Math.}\ {\bf 121} \#1 (1995), 211--233.

[GS2]  Garcia, A., Stichtenoth, H.: Asymptotically good towers of
function fields over finite fields.  {\em C. R. Acad.\ Sci.\ Paris~I}\/
{\bf 322} (1996), 1067--1070.

[I] Ihara, Y.: Some remarks on the number of rational points of
algebraic curves over finite fields.  {\em J. Fac.\ Sci.\ Tokyo}
{\bf 28} (1981), 721--724.

[S1] Sol\'e, P.: $D_4, E_6, E_8$ and the AGM.  Pages 448--455 in
{\em Applied Algebra, Algebraic Algorithms, and Error-Correcting Codes:
11th International Symposium (AAECC-11), Paris, France, July 1995}
(G.~Cohen, M.~Giusti, T.~Mora, eds.; Springer LNCS {\bf 948}).

[S2] Sol\'e, P.: Towers of Function Fields and Iterated Means.
Preprint, 1997.

[T] Takeuchi, K.: Commensurability classes of arithmetic triangle
groups.  {\em J.~Fac.\ Sci.\ Univ.\ Tokyo} {\bf 24} (1977), 201--212.

[TVZ] Tsfasman, M.A., \Vlad, S.G., Zink, T.:
Modular curves, Shimura curves and Goppa codes better than the 
Varshamov-Gilbert bound.  {\em Math.\ Nachr.}\ {\bf 109} (1982),
21--28.

\end{document}